\documentclass[]{amsart}

\usepackage[utf8]{inputenc}
\usepackage{amssymb}
\usepackage{amsmath}
\usepackage{amstext}
\usepackage{amsthm}
\usepackage{enumitem}

\usepackage{graphicx}
\usepackage{tabularx}

\usepackage[colorlinks=true]{hyperref} %

\usepackage{subcaption}

\usepackage[citestyle=alphabetic,bibstyle=alphabetic,backend=bibtex]{biblatex}
\usepackage{tikz}
\usepackage{tikz-cd}
\usetikzlibrary{decorations.pathreplacing}
\newcommand{\tikzmark}[1]{\tikz[overlay,remember picture,baseline=(#1.base)]
	\node (#1) {\strut};}

\usepackage[linesnumbered,ruled,vlined]{algorithm2e}

\SetCommentSty{mycommfont}

\newcommand{\tr}{\mathop{\mathrm{tr}}}

\newcommand{\spec}{\mathop{\mathrm{Spec}}}

\newcommand{\nil}{\mathop{\mathrm{nil}}}
\newcommand{\ideal}[1]{\langle#1\rangle}
\newcommand{\bs}[1]{\boldsymbol{#1}}

\newcommand{\lc}{{\mathtt{LC}}}

\newcommand{\sres}{{\mathtt{sRes}}}
\newcommand{\etale}{{\'{e}tale }}

\newtheorem{theorem}{Theorem}[section]
\newtheorem{result}{Main Result}

\newtheorem{lemma}[theorem]{Lemma}
\newtheorem{proposition}[theorem]{Proposition}

\theoremstyle{definition}
\newtheorem{definition}[theorem]{Definition}
\newtheorem{example}{Example}

\newtheorem*{problem*}{Problem}
\theoremstyle{remark}
\newtheorem{remark}{\indent Remark}
\newtheorem*{remark*}{Remark}

\title[Morphisms inducing Covering Maps over Real Closed Field]{What Kind of Morphisms Induces Covering Maps over a Real Closed Field?}
\author{Rizeng Chen}

\address{School of Mathematical Sciences, Peking University, 100871, Beijing, China}
\email{xiaxueqaq@stu.pku.edu.cn}

\keywords{Covering Map, \'Etale  Morphism, Geometric Fiber, Real Algebraic Geometry}

\subjclass[2020]{Primary 14P10; Secondary 14A10, 14Q30}
\thanks{This work is supported by National Key R\&D Program of China (No. 2022YFA1005102, No. 2024YFA1014003).}

\begin{document}

\maketitle

\begin{abstract}	
	In this article, we show that a flat morphism of $k$-varieties ($\mathop{\mathrm{char}} k=0$) with locally constant geometric fibers becomes finite \'etale after reduction. When $k$ is a real closed field, we prove that such a morphism induces a covering map on the rational points. We further give a triviality result different from Hardt's and a new interpretation of the construction of cylindrical algebraic decomposition as applications.
\end{abstract}

\section{Introduction}
\subsection{Problem Statement}
Covering map is probably one of the most important objects in topology. It is very useful in studying the fundamental group of a topological space and it also plays an important role in the theory of Riemann surfaces and complex manifolds. It is then natural to ask: what is the algebro-geometric analogue of covering maps? 

Of course, Zariski topology is too coarse for such a problem. So we shall restrict our attention to a smaller category of schemes with nicer topology. Fortunately, a real closed field $R$ have a stronger topology, namely the Euclidean topology. For an $R$-variety $X$, $X(R)$ can be naturally endowed with the Euclidean topology. Then we may ask:

\begin{problem*}\label{prob-main}
	In the category of $R$-varieties, what kind of morphisms induce a covering map on the rational points (equipped with Euclidean topology)?
\end{problem*}

The motivation behind this problem actually comes not only from pure mathematics, but also from other scientific fields. For example, the stability of a chemical reaction network or a biological system can be modeled by a polynomial system with parameters \cite{gatermann2001counting,craciun2009toric,wang2005stability}, and one wants to study what kind of parameters leads to a certain number of distinct (complex/real/positive) solutions. Notice that there is a canonical projection from the zero locus of the polynomial system to the space of parameters. So if the projection induces a covering map, then it is easy to tell the number of distinct solutions.

It seems that the best candidate answer should be finite \'etale morphisms. For example, when $R=\mathbb{R}$, a finite \'etale morphism of non-singular $\mathbb{R}$-varieties induces a covering map on the $\mathbb{R}$-rational points equipped with the Euclidean topology (Ehresmann's lemma). At the same time, for a finite \'etale surjective morphism $\varphi:X\to S$ with $S$ connected, there is a finite, locally free and surjective morphism $\psi:Y\to S$ such that $X\times_S Y$ is a trivial cover over $Y$ \cite[Proposition 5.2.9]{szamuely2009galois}. Grothendieck used finite \'etale surjective morphisms as the building block for the development of algebraic fundamental group \cite{grothendieck1971revetements}. %

However, finite \'etale map alone does not give a satisfactory answer.
\begin{enumerate}
	\item \'Etaleness might be stronger than needed. Consider a toy example, the double line $\spec R[x,y]/\ideal{y^2}$ over $\spec R[x]$. This is obviously a covering map since the underlying sets are homeomorphic. But it is not an \'etale map because the scheme-theoretic fibers are not smooth. This phenomenon is very common when studying parametric polynomial systems from the real world. 
	\item Less is known about the map on the rational points. Certainly, the Euclidean topology of $R$-points must be taken into consideration. One needs to develop specific methods to show that the induced map is a covering map.
\end{enumerate}

\subsection{Our Contribution}
In this paper, we will give an answer to Problem \ref{prob-main} addressing these two issues.

First, we propose a new family of morphisms to replace finite \'etale morphisms, namely the q-\'etale morphism.
\begin{definition}
	A morphism $f:X\to Y$ between $k$-varieties $X,Y$ is said to be \textbf{q-\'etale} if 
	\begin{itemize}
		\item $f$ is flat.
		\item for every $y\in Y$, there is an open neighborhood $U$ of $y$ such that the number of geometric points in $f^{-1}(y)$ is a finite constant on $U$.
	\end{itemize}
\end{definition}
The latter condition can be equivalently stated as ``the geometric fiber cardinality $n(y):y\mapsto\# X_{\overline{y}}$ is a locally constant function from $Y$ to $\mathbb{Z}$''. It is an analogue of the notion of unramified morphism: once the number of geometric points in the fiber is fixed, the fiber would not ``split''. If we replace it with $f$ being unramified, then $f$ is \'etale in the usual sense. Therefore the new definition is a plausible one. Such kind of morphisms arises naturally in computational algebraic geometry, where one studies solutions of a polynomial system with parameters.

It might be not easy to immediately see that q-\'etale morphisms have a deep connection to finite \'etale morphisms. The following result proved in the paper shows that q-\'etale morphisms possess the same topological properties as finite \'etale morphisms do over a field of characteristic 0, thus addressing the first issue.
\begin{result}[Theorem \ref{thm-reduced-g-etale-is-etale}]
	Let $k$ be a field of characteristic $0$. If $\pi:X\to Y$ is a q-\'etale $k$-morphism of $k$-varieties, then $\pi_{\mathrm{red}}:X_{\mathrm{red}}\to Y_{\mathrm{red}}$ is a finite \'etale morphism.
\end{result}
So for a flat morphism of $k$-varieties ($\mathop{\mathrm{char}}k=0$), if the cardinality of geometric fiber is a constant, then the ramification is purely due to the non-reduced structures of $X,Y$ and can be resolved by reduction. Moreover, the flatness is also preserved (notice that the reduction of a flat map is not necessarily flat). Proving that $\pi_{\mathrm{red}}$ is flat is actually the main difficulty in establishing Theorem \ref{thm-reduced-g-etale-is-etale}. Our strategy is to use a d\'evissage argument to show that $X_{\mathrm{red}}$ is locally isomorphic to the spectrum of a finite standard \'etale algebra over an affine open of $Y_{\mathrm{red}}$. %

Now focus on $R$-varieties. We show that q-\'etale morphisms do induce covering maps on $R$-rational points, thus addressing the second issue and answering our question.
\begin{result}[Theorem \ref{thm-covering-map}]
	Suppose $\pi: X\to Y$ is a q-\'etale $R$-morphism of $R$-varieties. Then $\pi_R:X(R)\to Y(R)$ is a covering map in the Euclidean topology.
\end{result}

This yields a triviality result similar to Hardt's \cite{hardt1980semi}, and it also explains the core mechanism of another classical construction in real algebraic geometry, the cylindrical algebraic decomposition \cite{collins1975quantifier}. We will discuss about this in more details in Section \ref{sect-app}. It is worth mentioning that while the preconditions of this result are purely geometric (flatness and locally constant geometric fiber), independent of the order structure on $R$, the conclusion is a semi-algebraic property (local existence of semi-algebraic sections of $\pi_R$). 

\subsection{Related Works}
	There are several related works concerning covering maps over a real closed field. Delfs and Manfred discussed covering maps in the category of locally semi-algebraic spaces \cite[Section 5]{delfs1984introduction}, generalizing the notion of semi-algebraic spaces to study the topology of semi-algebraic sets, especially their homology and homotopy \cite{delfs1981semialgebraic,knebusch1981semialgebraic,delfs1985locally,knebusch1989weakly}. Schwartz characterized a covering map $f$ of locally semi-algebraic spaces by the associated morphism of the corresponding real closed spaces being finite and flat \cite{schwartz1988open}. The theory of real closed spaces is introduced by Schwartz as a generalization of locally semi-algebraic spaces \cite{schwartz1989basic}. Baro, Fernando and Gamboa studied the relationship between branched coverings of semi-algebraic sets and its induced map on the spectrum of the ring of continuous semi-algebraic functions in \cite{baro2022spectral}. A common feature of these works is the use of (various generalizations) of semi-algebraic space. Hence, the sheaf equipped is the sheaf of continuous semi-algebraic functions. As our purpose is to study the properties of (scheme-theoretic) real varieties morphisms that induce covering maps, we will be using the sheaf of regular functions instead. %

	In \cite{bernard2024algebraic}, Bernard et al.\ studied the algebraic condition for a morphism between reduced varieties over $R[\sqrt{-1}]$ to be a homeomorphism for the Zariski topology or the Euclidean topology. This is similar to our interests, as both papers investigate what kind of variety morphisms have specific topological properties (being a homeomorphism or covering, respectively).
\section{Preliminary}

\subsection{Notations and Conventions}
\begin{itemize}
	\item In this paper, a \textbf{variety} is a separated, of finite type scheme over a field $k$.
	\item A \textbf{geometric point} $\overline{x}$ of a scheme $X$ is a morphism $\overline{x}: \spec \Omega\to X$, where $\Omega$ is an algebraically closed field. We say $\overline{x}$ lies over $x$ to indicate that $x\in X$ is the image of $\overline{x}$. 
	\item Given a morphism of schemes $f:X\to Y$ and a geometric point $\overline{y}:\spec \Omega\to Y$, the \textbf{geometric fiber} $X_{\overline{y}}$ is the fiber product $X\times_Y \spec \Omega$. %
	\item Suppose $f:X\to Y$ is a $k$-morphism of $k$-varieties and $K/k$ is a field extension, then $f_K:X(K)\to Y(K)$ is the map on the $K$-points.
	\item We fix $R$ to be a real closed field and $C=R[\sqrt{-1}]$ is the algebraic closure of $R$.
	\item By a \textbf{covering map}, we mean a continuous map $p:E\to B$ such that every $b\in B$ has an open neighborhood $U$ whose preimage is a disjoint (possibly empty) union of opens homeomorphic to $U$. We do not assume $p$ to be surjective here because $E=X(k)$ sometimes can be empty.
\end{itemize}

\subsection{Ingredients from Computational Algebraic Geometry}
With the emergence of computer, researchers became more and more interested in designing effective and constructive methods in algebraic geometry. These methods can be turned into algorithms to solve many problems, e.g., solving polynomial systems. This gave birth to the field of computational algebraic geometry, which has been a great source of ideas for decades. Here, we collect a few results that will be used later.
\subsubsection{Subresultants}

Subresultant is a classical tool in computational algebraic geometry. It characterizes the degree of the greatest common divisor of two univariate polynomials. 
\begin{definition}
	Let $A$ be a commutative ring with identity. Given two polynomials $f=\sum_{i=0}^{p} a_i x^i, g=\sum_{i=0}^{q} b_i x^i\in A[x]$ $(\deg f=p>q=\deg g)$. The \textbf{$j$-th Sylvester matrix} $\mathrm{Syl}_j(f,g)$ is the following matrix 

$$\left(\begin{array}{cccccccc}
	a_p     & 0   	  & \cdots & 0 & b_q 	 & 0   & \cdots & 0 \\
	a_{p-1} & a_p     & \cdots & 0 & b_{q-1} & b_q & \cdots & 0 \\
	a_{p-2} & a_{p-1} & \ddots & 0 & b_{q-2} & b_{q-1} & \ddots & 0 \\
	\vdots  &\vdots   & \ddots &a_p&\vdots   &\vdots& \ddots &b_q \\
	a_0     &a_{1}    & \cdots &\vdots &b_0& b_{1}& \cdots &\vdots \\
	0       &a_{0}    & \ddots &\vdots &0  & b_{0}  & \ddots &\vdots \\
	\vdots  &\vdots   & \ddots &a_{1} &\vdots&\vdots& \ddots &b_{1} \\
	\tikzmark{lower1L}0 & 0 & \cdots & a_0\tikzmark{lower1R} & \tikzmark{lower2L}0 & 0 & \cdots & b_0\tikzmark{lower2R}
\end{array}\right)
\begin{tikzpicture}[overlay, remember picture,decoration={brace,amplitude=5pt}]
	\draw[decorate,thick] (lower1R.south) -- (lower1L.south)
	node [midway,below=5pt] {$q-j$ columns};
	\draw[decorate,thick] (lower2R.south) -- (lower2L.south)
	node [midway,below=5pt] {$p-j$ columns};
\end{tikzpicture}
$$	
\vspace*{8mm}\\
representing the $A$-linear map $$\varphi_j:\begin{array}{rcl} A[x]_{q-j}\oplus A[x]_{p-j} & \to & A[x]_{p+q-j}
	\\ (u,v)&\mapsto& uf+vg.
\end{array}$$
between free modules in the canonical monomial bases.

Clearly this matrix has $(p+q-2j)$ columns and $(p+q-j)$ rows.

The \textbf{$j$-th subresultant polynomial} $\mathtt{sResP}_{x,j}(f,g)$ is the polynomial $\sum_{i=0}^{j}d_{j,i}{x}^i$ where $d_{j,i}$ is the $(p+q-2j)\times (p+q-2j)$ minor of $\mathrm{Syl}_{j}(f,g)$ extracted on the rows $1,\ldots,p+q-2j-1, p+q-j-i$. 

The leading coefficient $d_{j,j}$ of $\mathtt{sResP}_{x,j}(f,g)$, which is the minor of $\mathrm{Syl}_{j}(f,g)$ from the first $p+q-2j$ rows, is called the \textbf{$j$-th subresultant}, denoted by $\sres_{x,j}(f,g)$.
Especially, the $0$-th subresultant $\sres_{x,0}(f,g)$ is the resultant of $f$ and $g$.
\end{definition}

From the determinantal definition, it is immediate that subresultants and subresultant polynomials have the specialization property, provided that the leading coefficients are not mapped to zero.

\begin{lemma}
	Let $\varphi:A\to B$ be a ring homomorphism and $f,g\in A[x]$. Suppose that $\lc_x(f),\lc_x(g)\notin \ker \varphi$, then $$\varphi(\sres_{x,j}(f,g))=\sres_{x,j}(\varphi(f),\varphi(g))$$  and  $$\varphi(\mathtt{sResP}_{x,j}(f,g))=\mathtt{sResP}_{x,j}(\varphi(f),\varphi(g)).$$
\end{lemma}

The following classical result relates subresultants to the greatest common divisor of two univariate polynomials over a field. A proof may be found in \cite[Proposition 4.25, Theorem 8.56]{basu2008algorithms}.
\begin{theorem}
	Now suppose $A=k$ is a field. 
	The greatest common divisor of $f$ and $g$ in $k[x]$ is of degree $j$, if and only if
	$$\mathtt{sRes}_{x,0}(f,g)=\cdots=\mathtt{sRes}_{x,j-1}(f,g)=0, \mathtt{sRes}_{x,j}(f,g)\ne 0.$$
	
	Moreover, the $j$-th subresultant polynomial $\mathtt{sResP}_{x,j}$ coincides with the greatest common divisor of $f$ and $g$ in this case.
\end{theorem}

\subsubsection{Separating Elements}

Suppose $k$ is a field of characteristic 0. Let $A$ be a reduced $k$-algebra of finite type, and let $B=A[x_1,\ldots,x_n]/I$ be an $A$-algebra that is also a free $A$-module of finite rank $d$. Every element $\sigma\in B$ defines an $A$-linear multiplication map $L_{\sigma}:\begin{array}{rcl}
	B&\to &B \\ b&\mapsto & \sigma b
\end{array}$. Denote the characteristic polynomial $\det(\lambda I-L_{\sigma})\in A[\lambda]$ by $\chi_\sigma$. %

\begin{definition}
	An element $\sigma\in B$ is said to be \textbf{separating} over $p\in \spec A$, if $$\# (\spec B)_{\overline{p}}=d-\deg_\lambda\gcd(\chi_\sigma(p),\chi_\sigma'(p)).$$ This is equivalent to say the degree of square-free part of $\chi_\sigma(p)$ is equal to cardinality of geometric fiber over $p$. If $\sigma$ is separating over every $p\in \spec A$, then we say $\sigma$ is a \textbf{separating element} of $B$.
\end{definition}

At the first glance, it may be hard to see why such a condition is said to be separating. Stickelberger's Eigenvalue Theorem \cite[Theorem 4.2.7]{cox2005using} justifies this.
\begin{theorem}[Stickelberger]\label{thm-stickelberger}
	Let $k$ be an algebraically closed field.
	If $J$ is a zero-dimensional ideal in $k[x_1,\ldots,x_n]$ and $\sigma\in k[x_1,\ldots,x_n]$, then 
	$$\det(\lambda I-L_\sigma)=\prod_{p\in V(J)} (\lambda - \sigma(p))^{\mu(p)},$$
	where $\mu(p)=\dim_k (k[x_1,\ldots,x_n]/J)_p$ is the multiplicity of $p$ in $V(J)$ and $L_\sigma$ is the multiplication map by $\sigma $ in $k[x_1,\ldots,x_n]/J$. Moreover we have 
	$$\tr L_\sigma=\sum_{p\in V(J)}\mu(p)\sigma(p)\text{ and } \det L_\sigma=\prod_{p\in V(J)} \sigma(p)^{\mu(p)}.$$
\end{theorem}

By the Stickelberger's Eigenvalue Theorem, the degree of the square-free part of $\chi_\sigma(p)$ is less than or equal to the number of points in the geometric fiber $(\spec B)_{\overline{p}}$, and the equality holds if and only if $\sigma(p)$ takes distinct values on each point of the geometric fiber. In this case, $\sigma$ separates different points in $(\spec B)_{\overline{p}}$.

The following well-known lemma shows that there is always a separating element for a zero-dimensional $k$-varieties when $\mathop{\mathrm{char}}k=0$. Readers may refer to \cite[Lemma 4.90]{basu2008algorithms} for a proof.

\begin{lemma}\label{lem-enough-separating-elements}
	Suppose $k$ is a field of characteristic 0. Let $A=k[x_1,\ldots,x_n]/J$ be a zero-dimensional $k$-algebra such that $\dim_k A=d$. Among the elements of 
	$$\left\{ x_1+ix_2+\cdots+i^{n-1}x_n \right\}_{0\leq i\leq (n-1)\binom{d}{2}},$$
	at least one of them is a separating element for $A$.
\end{lemma}

\subsubsection{Rouillier's Lemma}

The following lemma due to Fabrice Rouillier appeared in his proof of existence of Rational Univariate Representation \cite{rouillier1999solving}. We include the proof here because it is succinct.
\begin{lemma}[Rouillier]\label{lem-rouillier}
	Let $k$ be an algebraically closed field of characteristic 0.
	Suppose $A=k[x_1,\ldots,x_n]/J$ is zero-dimensional, $\dim_k A=d$ and $\sigma\in k[x_1,\ldots,x_n]$ is a separating element of $A$. Let $u=\sum_{i=0}^{d}u_i \lambda^i$ be the monic square-free part of $\chi_\sigma=\det(\lambda I-L_\sigma)$. For any $f\in k[x_1,\ldots,x_n]$, we have 
	$$\sum_{l=0}^{d-1}\sum_{i=0}^{d-l-1}\tr L_{f\sigma^i}\cdot u_{l+i+1}\cdot\lambda^l = 
	\sum_{\alpha\in V(J)} \mu(\alpha)f(\alpha)\prod_{\beta\in V(J),\beta\ne\alpha}(\lambda-\sigma(\beta)).$$
\end{lemma}
\begin{proof}
	Notice that $u=\prod_{\beta \in V(J)}(\lambda-\sigma(\beta))$. So we have the following computation in $k((\lambda^{-1}))$ (the field of formal Laurent series in $\lambda^{-1}$):
	\begin{align*}
		~&\sum_{\alpha\in V(J)} \mu(\alpha)f(\alpha)\prod_{\beta\in V(J),\beta\ne\alpha}(\lambda-\sigma(\beta))\\
		=&\left(\prod_{\beta \in V(J)}\left(\lambda-\sigma(\beta)\right)\right)\left(\sum_{\alpha \in V(J)} \frac{\mu(\alpha)f(\alpha)}{\lambda-\sigma(\alpha)}\right) &&(\text{Common factor})\\
		=&\left(\prod_{\beta \in V(J)}\left(\lambda-\sigma(\beta)\right)\right)\left(\sum_{\alpha \in V(J)}\sum_{i=0}^{+\infty} \frac{\mu(\alpha)f(\alpha)\sigma(\alpha)^i}{\lambda^{i+1}}\right) && \text{(Series expansion)}\\
		=&\left(\sum_{j=0}^{d} u_j\lambda^j\right)\left(\sum_{i=0}^{+\infty}\lambda^{-i-1}\sum_{\alpha \in V(J)} \mu(\alpha)f(\alpha)\sigma(\alpha)^i\right) && (\sigma\text{ separating})\\
		=&\left(\sum_{j=0}^{d} u_j\lambda^j\right)\left(\sum_{i=0}^{+\infty}\lambda^{-i-1}\tr L_{f\sigma^i}\right)&& \text{(Stickelberger)} \\
		=&\sum_{l=0}^{d-1}\sum_{i=0}^{d-l-1} u_{l+i+1}\cdot \tr L_{f\sigma^i}\cdot \lambda^l + \sum_{l<0}\sum_{i=-l-1}^{d-l-1} u_{l+i+1}\cdot \tr L_{f\sigma^i}\cdot \lambda^l.		
	\end{align*}
	Now it suffices to show that the second term is zero. Observe that when $l<0$:
	$$\begin{aligned}[t]\sum_{i=-l-1}^{d-l-1}  u_{l+i+1}\tr L_{f\sigma^i}=\sum_{i=0}^{d} u_{i}\tr L_{f\sigma^{i-l-1}}=\sum_{i=0}^{d}u_i\sum_{\alpha\in V(J)}\mu(\alpha)f(\alpha)\sigma^{i-l-1}(\alpha)\\
		=\sum_{\alpha\in V(J)}\left(\sum_{i=0}^d u_i \sigma^i(\alpha)\right)\mu(\alpha)f(\alpha)\sigma^{-l-1}(\alpha)
		=\sum_{\alpha\in V(J)}u(\sigma(\alpha))\mu(\alpha)f(\alpha)\sigma^{-l-1}(\alpha).
		\end{aligned}$$
	Since $u$ nullifies $\sigma(\alpha)$ by Stickelberger's Eigenvalue Theorem, the above sum is zero and our assertion is proved.
\end{proof}

\section{Reduction of q-\'{e}tale Morphism}
	In this section, we will prove our first main result: a q-\'etale morphism between $k$-varieties ($\mathop{\mathrm{char}}k=0$) becomes finite \'etale after reduction.

	The finiteness can be derived from the following result \cite[Proposition 15.5.9]{grothendieck1966elements}.
	\begin{proposition}[Grothendieck]\label{prop-grothendieck}
		Let $f:X\to Y$ be a flat morphism locally of finite presentation. For each $y\in Y$, let $n(y)$ be the number of geometric connected components of $f^{-1}(y)$.
		\begin{enumerate}[label*=(\roman*)]
			\item If $f$ is quasi-finite and separated, then the function $y\mapsto n(y)$ (which is the number of geometric points in $f^{-1}(y)$) is lower semi-continuous on $Y$. If it is a constant in a neighborhood of $y_0$, then $f$ is proper on a neighborhood of $y_0$.
			\item If $f$ is proper, then function $y\mapsto n(y)$ is lower semi-continuous. Suppose moreover $f^{-1}(y_0)$ is geometrically reduced over $\kappa_{y_0}$, then $n(y)$ is a constant on a neighborhood of $y_0$.
		\end{enumerate}
	\end{proposition}
	
	\begin{proposition}\label{prop-qff-implies-ff} %
 		A q-\'etale $k$-morphism of $k$-varieties is finite.	
	\end{proposition}
	\begin{proof}
		Let $f:X\to Y$ be such a morphism.
		By Proposition \ref{prop-grothendieck}, for any $y\in Y$, there is a neighborhood $U$ of $y$ such that $f^{-1}(U)\to U$ is proper. Notice that a quasi-finite proper morphism is finite \cite[Th\'eor\`eme 8.11.1]{grothendieck1966elements}. Therefore $f^{-1}(U)\to U$ is finite, and so is $f$.
	\end{proof}
	
	Conversely, if $\pi:X\to Y$ is a finite \'etale $k$-morphism of $k$-varieties, then $\pi$ is q-\'etale. This is an easy corollary of \cite[Proposition 18.2.8]{grothendieck1967elements}, which is shown below:%
	\begin{proposition}[Grothendieck]
		Let $f:X\to Y$ be a separated, of finite type, \etale morphism. For each $y\in Y$, let $n(y)$ be the number of geometric points of $f^{-1}(y)$. Then the function $n(y)$ is lower semi-continuous on $Y$. For it to be continuous near $y$ (i.e.\ to be a constant in a neighborhood of $y$), it is necessary and sufficient to have an open neighborhood $U$ of $y$ such that the restriction $f^{-1}(U)\to U$ of $f$ is a finite (\'etale) morphism.
	\end{proposition}
	
	Now all the ingredients for the first main result are ready.
\begin{theorem}\label{thm-reduced-g-etale-is-etale}
	Let $k$ be a field of characteristic 0. Suppose $X,Y$ are $k$-varieties and $\varphi:X\to Y$ is a q-\'etale $k$-morphism, then the reduced map $\varphi_{\mathrm{red}}: X_{\mathrm{red}} \to Y_{\mathrm{red}}$ is (finite) \'etale.
	
	Moreover, the reduced map $\varphi_{\mathrm{red}}$ is locally of the form $\spec A[\lambda]/\ideal{u}\to \spec A$, where $u\in A[\lambda]$ is a monic polynomial and $u'$ is invertible modulo $u$.
\end{theorem}
\begin{proof}The proof is a standard dévissage argument. Note that the finiteness is already shown in Proposition \ref{prop-qff-implies-ff}. \medskip
	
	\textit{Step 1. Reduce to the affine case.} The question is local, we may assume $X=\spec B$, $Y=\spec A$ are affine varieties and $B=A[x_1,\ldots,x_n]/I$ is a free $A$-module of finite rank $r$ because $\varphi$ is finite flat by Proposition \ref{prop-qff-implies-ff}. Replacing $A$, $B$ with $A_{\mathrm{red}}=A/\nil A$ and $B\otimes_A A_{\mathrm{red}}$ respectively, we may even assume that $A$ is reduced. \medskip
	
	\textit{Step 2. Find a separating element.} Fix arbitrary $y\in Y$. By the existence of separating element (Lemma \ref{lem-enough-separating-elements}), there is $\sigma \in B$ taking different values on geometric points of $\varphi^{-1}(y)$. That is, $$r-\deg_\lambda\gcd(\chi_\sigma(y;\lambda),\chi_\sigma'(y;\lambda))=\# X_{\overline{y}},$$ where $\chi_\sigma=\det(\lambda I-L_\sigma)$. Let $s=r-\# X_{\overline{y}}$.
	
	Now let $\Delta_i=\sres_{\lambda,i}(\chi_\sigma,\chi_\sigma')\in A$ for $i=0,\ldots,r-1$. Then by the specialization property of subresultants, $\Delta_0,\cdots,\Delta_{s-1}$ vanish everywhere on $Y$ but $\Delta_s(y)\ne 0$. Since $A$ is reduced, we have that $\Delta_0=\cdots=\Delta_{s-1}=0$ and $D(\Delta_s)$ is a non-empty affine open neighborhood of $y\in Y$. Replacing $A$, $B$ with $A_{\Delta_s}$ and $B_{\Delta_s}$, we may assume that $\sigma \in B$ is a separating element. \medskip
	
	\textit{Step 3. Intermediate closed subscheme of $\mathbb{A}_Y^1$.} Let $Z=\spec A[\lambda]/\ideal{\chi_\sigma}$. Since $\chi_\sigma(L_\sigma)=0$ by Cayley-Hamilton Theorem, we see that the ring map $\begin{array}{rcl}A[\lambda]&\to& B\\\lambda &\mapsto&\sigma \end{array}$ factors through $A[\lambda]/\ideal{\chi_\sigma}$. Therefore there are $k$-morphisms of $k$-varieties making the following diagram commute.
	\begin{center}
		\begin{tikzcd}
			X \arrow[dr, "\varphi"]\arrow[d,swap,"\tilde{\varphi}"]&  \\
			Z \arrow[r,swap,"\pi"] & Y\\
		\end{tikzcd}
	\end{center}
	We claim that $\tilde{\varphi}_{\overline{k}}: X(\overline{k})\to Z(\overline{k})$ is bijective. To see this, we write down the coordinates explicitly $$\begin{array}{rcccl}X(\overline{k})&\to& Z(\overline{k}) &\to& Y(\overline{k})\\ (\bs{y},\bs{x})&\mapsto&(\bs{y},\lambda)=(\bs{y},\sigma(\bs{y},\bs{x}))&\mapsto&\bs{y}.\end{array}$$
	
 		On the one hand, the injectivity is obvious, because $\sigma$ is a separating element. On the other hand, let $(\bs{y},\lambda) \in Z(\overline{k})$, then $\lambda$ is an eigenvalue of $L_\sigma$ specialized at $\bs{y}$, so there is some $\bs{x}\in \varphi^{-1}_{\overline{k}}(\bs{y})$ such that $\sigma(\bs{y},\bs{x})=\lambda$ by Stickelberger's eigenvalue theorem. This shows the surjectivity. \medskip
 		
 	\textit{Step 4. The factorization of $\chi_\sigma$ and the reduced scheme associated to $Z$. } Now we will show that $\chi_\sigma$ is the product of two monic polynomials $u,f\in A[\lambda]$, where $u$ specializes to the square-free part of $\chi_\sigma(y;\lambda)$ and $f$ specializes to $\gcd(\chi_\sigma(y;\lambda),\chi_\sigma'(y;\lambda))$ for all specialization $y\in Y$. Moreover $Z_{\mathrm{red}}=\spec A[\lambda]/\ideal{u}$.
 	
	Consider the subresultants $\Delta_0,\ldots,\Delta_s$ (recall Step 2). Now $\Delta_0=\cdots=\Delta_{s-1}=0$ and $\Delta_{s}$ is a unit in $A$. In this case, the $s$-th subresultant polynomial $\mathtt{sResP}_{\lambda,s}(\chi_\sigma,\chi'_\sigma)\in A[\lambda]$, whose leading coefficient is $\Delta_{s}$, specializes to $\gcd(\chi_\sigma(p;\lambda),\chi_\sigma'(p;\lambda))$ for all $p\in \spec A$. 
 	
 	Let $f=\frac{1}{\Delta_{s}}\mathtt{sResP}_{\lambda,s}(\chi_\sigma,\chi'_\sigma)$. Then by performing polynomial division in $A[\lambda]$, %
 	there are $u,\gamma \in A[\lambda]$ ($\deg_\lambda \gamma<\deg_\lambda f$) such that $$\chi_\sigma=u\cdot f+\gamma.$$ But for every $p\in \spec A$, $f(p;\lambda)$ divides $\chi_\sigma(p;\lambda)$, implying that $\gamma \in p$. Hence, we have that $\gamma=0$ and $\chi_\sigma=uf$ is the desired factorization.
 	
	The assertion $Z_{\mathrm{red}}=\spec A[\lambda]/\ideal{u}$, or equivalently $\sqrt{\ideal{\chi_\sigma}}=\ideal{u}$, follows from these two observations. 
	\begin{itemize}
 		\item On the one hand, $u\in \sqrt{\ideal{\chi_\sigma}}$. Divide $u^r$ by $\chi_\sigma$ and denote the remainder by $\gamma$. For every specialization $p\in \spec A$, $\chi_\sigma(p;\lambda)|u^r(p;\lambda)$, so $\gamma(p;\lambda)=0$. Therefore $\gamma=0$ and $u^r\in \ideal{\chi_\sigma}$.
 		\item On the other hand, $u$ divides every element $g\in \sqrt{\ideal{\chi_\sigma}}$. Suppose $g^t\in\ideal{\chi_\sigma}$ and $g=uq+\gamma$ ($\deg_\lambda \gamma<\deg_\lambda u$).
		Then 
 		$$g^t=(uq+\gamma)^t=u(u^{t-1}q^t+\cdots+tq\gamma^{t-1})+\gamma^t\in \ideal{\chi_\sigma}\subseteq \ideal{u},$$
 		shows that $\gamma^t\in \ideal{u}$. If $\gamma$ is not zero then there is some $p\in \spec A$ such that $\gamma\notin p$. But $u(p;\lambda)$ is square-free, $u(p;\lambda)|\gamma^t(p;\lambda)$ forces $u(p;\lambda)|\gamma(p;\lambda)$.	
 		Then $\deg_\lambda \gamma\geq \deg_\lambda\gamma(p;\lambda)\geq \deg_\lambda u(p;\lambda)=\deg_\lambda u$, which is a contradiction.
 	\end{itemize}   \medskip
 	
 	\textit{Step 5. Construct an explicit isomorphism $X_{\mathrm{red}}\cong Z_{\mathrm{red}}$.} 	
 	Here, we will extend Rouillier's idea of rational univariate representation \cite{rouillier1999solving}, which gives an explicit bijection from solutions of a zero-dimensional multivariate system to roots of a univariate equation. Recall that $B=A[x_1,\ldots,x_n]/I$.
 	
	Let $g,g_1,\ldots,g_n\in A[\lambda]$ be the polynomials defined by $$g_i=\sum_{l=0}^{r-1}\sum_{j=0}^{r-l-1}\tr L_{x_i\sigma^j}\cdot u_{l+j+1}\cdot\lambda^l \text{ and } g=\sum_{l=0}^{r-1}\sum_{j=0}^{r-l-1}\tr L_{\sigma^j}\cdot u_{l+j+1}\cdot\lambda^l,$$ where $u=\sum_{i=0}^{r}u_i\lambda^i$ is the monic square-free part of $\chi_\sigma$.
	\begin{enumerate}[label*=(Step 5\alph*)]
		\item The $A$-algebra map $\begin{array}{rcl} A[x_1,\ldots,x_n] & \to & A[\lambda]/\ideal{\chi_\sigma} \\ x_i & \mapsto & g_ig^{-1}\end{array}$ is well-defined and gives a $Y$-morphism $\psi: Z \to \mathbb{A}_Y^n$. 
		
		To show the claim, we need to prove that $g$ is invertible modulo $\chi_\sigma$. For any $\overline{k}$-point $(\bs{y}_0,\lambda_0)$ of $Z$, it is the image of some $\overline{k}$-point $(\bs{y}_0,\bs{x}_0)$ of $X$. Now apply Lemma \ref{lem-rouillier} by setting $f=1$, we have 
		$$g(\bs{y}_0;\lambda)=\sum_{(\bs{y}_0,\bs{x})\in X(\overline{k})} \mu(\bs{y}_0,\bs{x})\prod_{\substack{(\bs{y}_0,\bs{x}')\in X(\overline{k}),\\ \bs{x}'\ne\bs{x}}}(\lambda-\sigma(\bs{y}_0,\bs{x}')),$$
		where $\mu(\bs{y}_0,\bs{x})$ is the multiplicity of $\bs{x}$ in the zero-dimensional fiber $X_{\bs{y}_0}$.
		
		Since $\lambda_0=\sigma(\bs{y}_0,x_0)$, we have $$g(\bs{y}_0,\lambda_0)=\mu(\bs{y}_0,\bs{x}_0)\prod_{\substack{(\bs{y}_0,\bs{x}')\in X(\overline{k}),\\ \bs{x}'\ne\bs{x}_0}}(\sigma(\bs{y}_0,\bs{x}_0)-\sigma(\bs{y}_0,\bs{x}'))\ne 0.$$
		So $g$ never vanishes on $Z(\overline{k})$, we conclude that $g$ is invertible modulo $\chi_\sigma$ by Nullstellensatz.
		\item Next, we claim that for any $\overline{k}$-point $(\bs{y}_0,\bs{x}_0)\in X(\overline{k})$, $g_i(\tilde{\varphi})/g(\tilde{\varphi})$ evaluates to $x_i$ for $1\leq i\leq n$. Let $\lambda_0=\sigma(\bs{y}_0,\bs{x}_0)$.

		To show the claim, we apply Lemma \ref{lem-rouillier} again by setting $f=x_i$. 
		$$g_i(\bs{y}_0;\lambda)=\sum_{(\bs{y}_0,\bs{x})\in X(\overline{k})} \mu(\bs{y}_0,\bs{x})x_i(\bs{y}_0,\bs{x})\prod_{\substack{(\bs{y}_0,\bs{x}')\in X(\overline{k}),\\ \bs{x}'\ne\bs{x}}}(\lambda-\sigma(\bs{y}_0,\bs{x}')),$$
		
		Therefore, $$g_i(\bs{y}_0,\lambda_0)=\mu(\bs{y}_0,\bs{x}_0)x_i(\bs{y}_0,\bs{x}_0)\prod_{\substack{(\bs{y}_0,\bs{x}')\in X(\overline{k}),\\ \bs{x}'\ne\bs{x}_0}}(\sigma(\bs{y}_0,\bs{x}_0)-\sigma(\bs{y}_0,\bs{x}')).$$
		
		So, we do have $\frac{g_i(\bs{y}_0,\sigma(\bs{y}_0,\bs{x}_0))}{g(\bs{y}_0,\sigma(\bs{y}_0,\bs{x}_0))}=x_i(\bs{y}_0,\bs{x}_0)$ for every geometric point $(\bs{y}_0,\bs{x}_0)\in X(\overline{k})$. Hence $g_i(\tilde{\varphi})/g(\tilde{\varphi})$ agrees with $x_i$ everywhere, implying that the induced morphism on reduced schemes $\psi_{\mathrm{red}}\circ \tilde{\varphi}_{\mathrm{red}}: X_{\mathrm{red}}\to Z_{\mathrm{red}}\to \mathbb{A}_Y^n$ is actually the closed immersion $\imath:X_{\mathrm{red}}\to \mathbb{A}_Y^n$. 
		
		\begin{center}
			\begin{tikzcd}
				X_{\mathrm{red}} \arrow[dr,hook', "\imath"]\arrow[d,swap,"\tilde{\varphi}_{\mathrm{red}}"]&  \\
				Z_{\mathrm{red}} \arrow[r,swap,"\psi_{\mathrm{red}}"] & \mathbb{A}_Y^n\\
			\end{tikzcd}
		\end{center}
		
		\item Now by shrinking the image, we have $\psi_{\mathrm{red}}\circ \tilde{\varphi}_{\mathrm{red}}=\mathrm{id}_{X_\mathrm{red}}$. Notice that $\tilde{\varphi}_{\mathrm{red}}$ is bijective on $\overline{k}$-points (because $\tilde{\varphi}$ is), therefore $\tilde{\varphi}_{\mathrm{red}}\circ \psi_{\mathrm{red}}$ must agree with $\mathrm{id}_{Z_{\mathrm{red}}}$ on $\overline{k}$-points too. Because two morphisms from a geometrically reduced scheme agreeing on geometric points are the same \cite[Exercise 5.17]{gortz2010algebraic}, we have $\tilde{\varphi}_{\mathrm{red}}\circ \psi_{\mathrm{red}}=\mathrm{id}_{Z_{\mathrm{red}}}$. Therefore $\tilde{\varphi}_{\mathrm{red}}$ is an isomorphism.
		
	\end{enumerate} \medskip
 	
 	\textit{Step 6. Complete the proof by showing $\pi_{\mathrm{red}}$ is \'etale.}
 	The flatness is clear, because $u\in A[\lambda]$ is a monic polynomial. For every $y\in Y$, $\# Z_{\overline{y}}=\# X_{\overline{y}}=\deg_\lambda u$, so the geometric fiber $(Z_{\mathrm{red}})_{\overline{y}}$ is zero-dimensional and regular. Therefore $\pi_{\mathrm{red}}:Z_{\mathrm{red}}\to Y$ is \'etale \cite[Theorem III.10.2]{hartshorne2013algebraic}. Because $\tilde{\varphi}_{\mathrm{red}}:X_{\mathrm{red}}\to Z_{\mathrm{red}}$ is an isomorphism, the proof is completed now.
\end{proof}

The last assertion of Theorem \ref{thm-reduced-g-etale-is-etale} is a slightly different version of the fact that \etale morphisms are locally standard \'etale.

\begin{example}
	Suppose $k$ is a field of characteristic 0. Consider $$X=\spec k[p,q,x]/\ideal{4p^3+27q^2,x^3+px+q}, Y=\spec k[p,q]/\ideal{4p^3+27q^2}$$ and the canonical projection $\pi:X\to Y$. Clearly $\pi$ is finite and flat. 
	
	\begin{figure}[hbtp]
		\centering
		\includegraphics[width=0.33\linewidth]{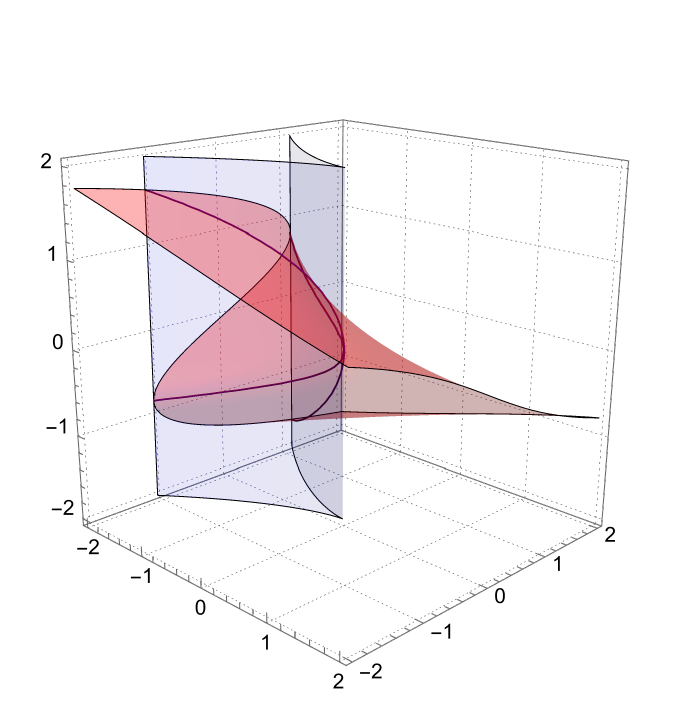}
		\caption{$X$ as the intersection of zeros of defining equations} 
		\label{fig:depressed-cubic-over-discriminant}
	\end{figure}
	
	Let $U=D(p)\subsetneq Y$, then $\pi|_{\pi^{-1}(U)}: X\times_Y U\to U$ is finite, flat and $\# X_{\overline{y}}=2$ for all $y\in U$. So by Theorem \ref{thm-reduced-g-etale-is-etale}, the reduced map $(X\times_Y U)_{\mathrm{red}}\to U_{\mathrm{red}}$ is finite \'etale.
	
	Notice that in general, even the flatness is not preserved under reduction. In this example, $\pi_{\mathrm{red}}: X_{\mathrm{red}}\to Y_{\mathrm{red}}$ is not flat! In fact, by a direct computation, $$\sqrt{\ideal{4p^3+27q^2,x^3+px+q}}=\ideal{q + p x + x^3, 4 p^2 - 9 q x + 6 p x^2,4p^3+27q^2},$$ so every fiber of $\pi_{\mathrm{red}}$ is of degree $2$ except the origin's, which is of degree $3$.
\end{example}

\begin{example}
	When $\mathop{\mathrm{char}} k=p>0$, the same conclusion does not hold anymore. Let $X=\spec \mathbb{F}_2[t][x]/\ideal{x^2+t}$ and $Y=\spec \mathbb{F}_2[t]$, then $X$ and $Y$ are both reduced, the canonical projection $\pi:X\to Y$ is flat and the cardinality of geometric fiber is always $1$. But $\pi$ is not smooth.
\end{example}	
	
\section{Semi-algebraic Covering Map}
Now we continue to prove our second main result.

\begin{theorem}\label{thm-covering-map}
	Let $X$, $Y$ be two $R$-varieties. Suppose $\pi: X\to Y$ is q-\etale $R$-morphism. Let $y\in Y(R)$. Then there exists $r\in \mathbb{N}$, semi-algebraically connected Euclidean neighborhood $U\subseteq Y(R)$ of $y$ and semi-algebraic continuous functions $\xi_1,\ldots,\xi_r:U\to X(R)$ such that
	\begin{enumerate}[label=(\roman*)]
		\item The number of real fibers over any $y\in U$, is always $r$.
		\item $\xi_i$ are sections of $\pi$, that is $\pi\circ\xi_i=\mathrm{id}_U$ for $i=1,\ldots,r$.
		\item For every $y\in U$: $i\ne j\implies \xi_i(y)\ne \xi_j(y)$.
	\end{enumerate}
	Therefore $\pi_R$ is a covering map in the Euclidean topology.
\end{theorem}
\begin{proof}
	We may replace $\pi:X\to Y$ with $\pi_{\mathrm{red}}:X_{\mathrm{red}}\to Y_{\mathrm{red}}$ and assume that $\pi$ is a finite \etale morphism by Theorem \ref{thm-reduced-g-etale-is-etale}. Then locally $\pi$ is of the form $\spec A[\lambda]/\ideal{u}\to \spec A$, where $u\in A[\lambda]$ is a monic polynomial and $u'$ is invertible modulo $u$. By the continuity of roots in the coefficients \cite[Theorem 5.12]{basu2008algorithms}, the number of distinct real roots of $u(\bs{y};\lambda)=0$ is a constant $r$ near $\bs{y}$, because the number of distinct complex roots (=the geometric fiber cardinality) is a constant. Intuitively, the continuity and the constancy of number of distinct roots guarantee that two different roots cannot collide and one multiple root cannot split when the coefficients vary. One may refer to \cite[Theorem 1]{collins1975quantifier} or \cite[Theorem 5.4]{chen2023geometric} for a detailed argument. %
	
	Define $$\xi_i:\begin{array}{rcl}U&\to &X(R)\\ \bs{y}&\mapsto& (\bs{y},\lambda_i),\end{array}$$
	where $\lambda_i$ is the $i$-th largest root of $u(\bs{y};\lambda)=0$ ($i=1,\ldots,r$), then $\xi_i$ are semi-algebraic functions by definition. Since the roots are continuous with respect to the polynomial coefficients, $\xi_i$ are continuous. The proof is completed. %
\end{proof}

\begin{remark}
	For readers familiar with locally semi-algebraic spaces, it can be immediately concluded that $\pi_R$ is a covering map after applying Theorem \ref{thm-reduced-g-etale-is-etale}, see \cite[Example 5.5]{delfs1984introduction}.
\end{remark}

\begin{example} \label{examples-of-covering-maps}
	Here are some examples fulfilling the q-\'etale condition. We see that they are indeed covering maps.
	\begin{figure}[htbp]
		\centering		
		\begin{subfigure}[t]{0.3\textwidth}
			\includegraphics[width=\linewidth]{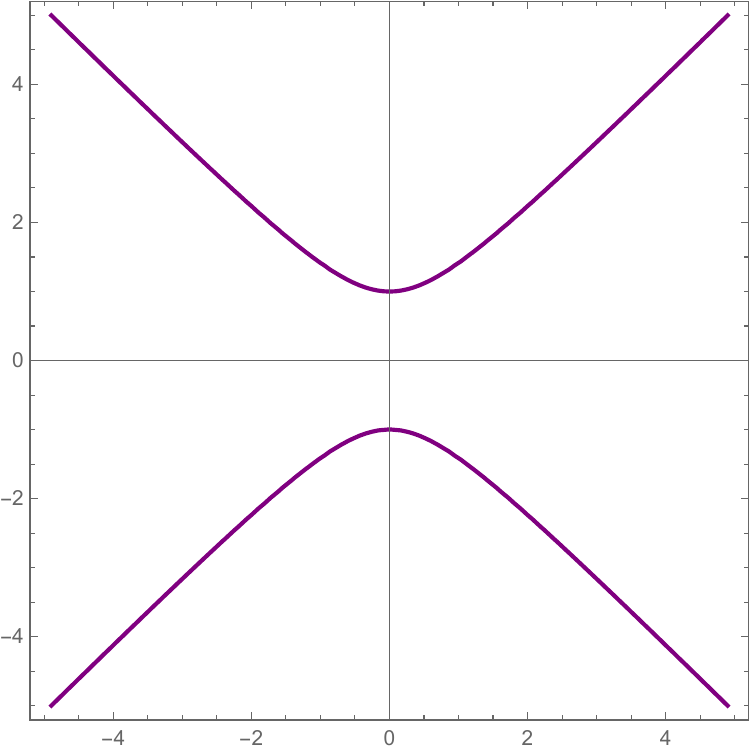}
			\caption{The curve defined by $y^2=x^2+1$} 
			\label{fig:y-squared-minus-x-squared-minus-one}
		\end{subfigure}
		~
		\begin{subfigure}[t]{0.3\textwidth}
			\includegraphics[width=\linewidth]{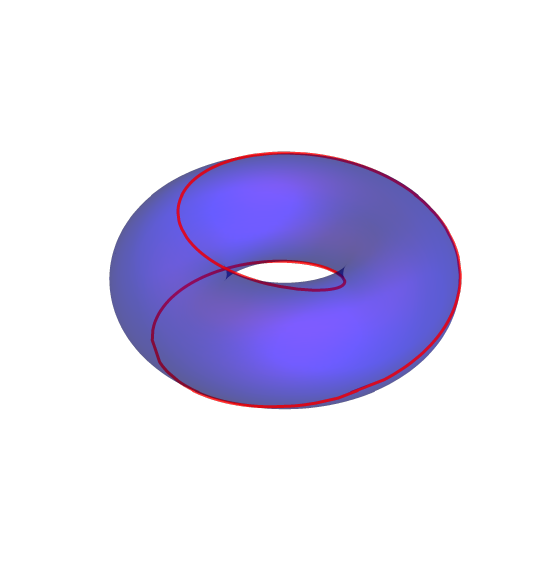}
			\caption{A curve on torus} 
			\label{fig:no-global-section}
		\end{subfigure}
		~
		\caption{Examples}
	\end{figure}	
	\begin{enumerate}[label=(\alph*)]
		\item Consider $X=\spec R[x,y]/\ideal{y^2-x^2-1}$ and $Y=\spec R[x]$. The projection is q-\'etale after removing $\ideal{x^2+1}\in \spec R[x]$, which is a pair of conjugate non-real points. Therefore $X(R)\to Y(R)$ is a covering map. See Figure \ref{fig:y-squared-minus-x-squared-minus-one}.
		\item Consider $\mathbb{R}$-varieties $X=\spec \mathbb{R}[x,y,z,w]/\ideal{x^2+y^2-1,z^2+w^2-1,x-2z^2+1,x+2 w^2-1,y-2 z w}$ and $Y=\spec \mathbb{R}[x,y]/\ideal{x^2+y^2-1}$. The canonical projection $\pi:X\to Y$ is flat with constant geometric fiber size $2$. The $\mathbb{R}$-rational points of $Y$ can be modeled by the unit circle $S^1$ and $X(\mathbb{R})$ can be identified with a curve parameterized by
		$\varphi\mapsto(\cos 2\varphi,\sin 2\varphi,\cos \varphi,\sin\varphi)$ in the torus $S^1\times S^1=\{(\cos \theta,\sin \theta,\cos \varphi,\sin \varphi)|\theta,\varphi\in\mathbb{R}\}$. So $\pi_\mathbb{R}$ is given by		
		$$(\cos 2\varphi,\sin 2\varphi,\cos \varphi,\sin\varphi) \mapsto (\cos \theta,\sin \theta,\cos \varphi,\sin \varphi)\mapsto(\cos \theta,\sin \theta).$$ 
		See Figure \ref{fig:no-global-section} for an embedding of $X(\mathbb{R})$ in $\mathbb{R}^3$ via the torus.
		
		In this example, one can only defined the semi-algebraic sections locally and $X(\mathbb{R})$ is not a disjoint union of copies of $Y(\mathbb{R})$. It is easy to see that the source is connected, but two different global sections will give two disjoint closed subsets.
	\end{enumerate}
\end{example}

The only assumption of Theorem \ref{thm-covering-map} is being q-\'etale, which is independent of the order structure of $R$, yet it is sufficient to conclude the local existence of semi-algebraic sections. Of course, such a simple condition is not necessary, as the map $x\mapsto x^3$ is a semi-algebraic homeomorphism but the geometric fiber cardinality changes when $x=0$. Also, one may add embedded points to destroy the flatness without changing the map on the topological level. Still, the below (non-)examples show that our result is quite sharp.

\begin{example}
	Theorem \ref{thm-covering-map} is no longer true without the q-\'etale hypothesis. We present some non-examples here.
	\begin{figure}[htbp]
		\centering		
		\begin{subfigure}[t]{0.3\textwidth}
			\includegraphics[width=\linewidth]{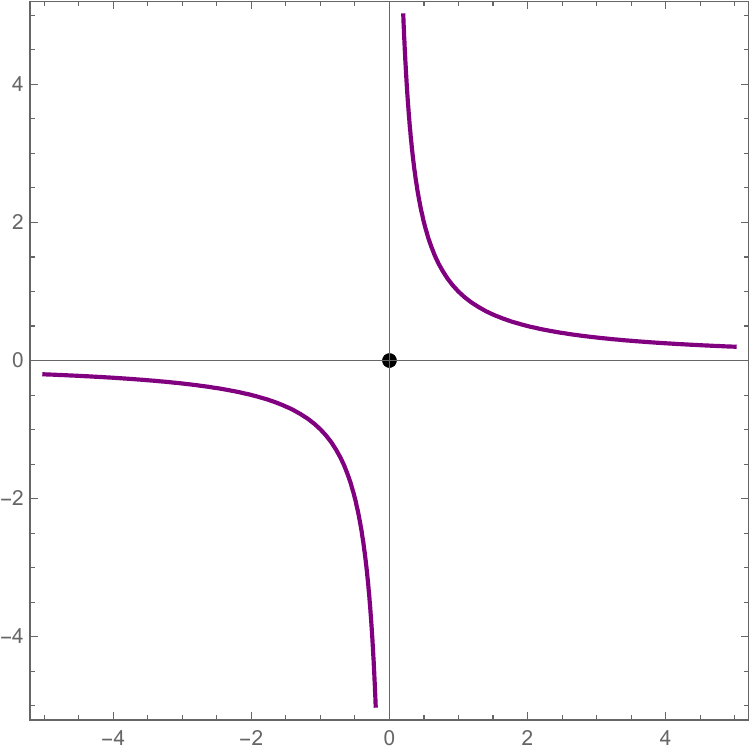}
			\caption{Without flatness} 
			\label{fig:non-flatness}
		\end{subfigure}
		~			
		\begin{subfigure}[t]{0.3\textwidth}
			\includegraphics[width=\linewidth]{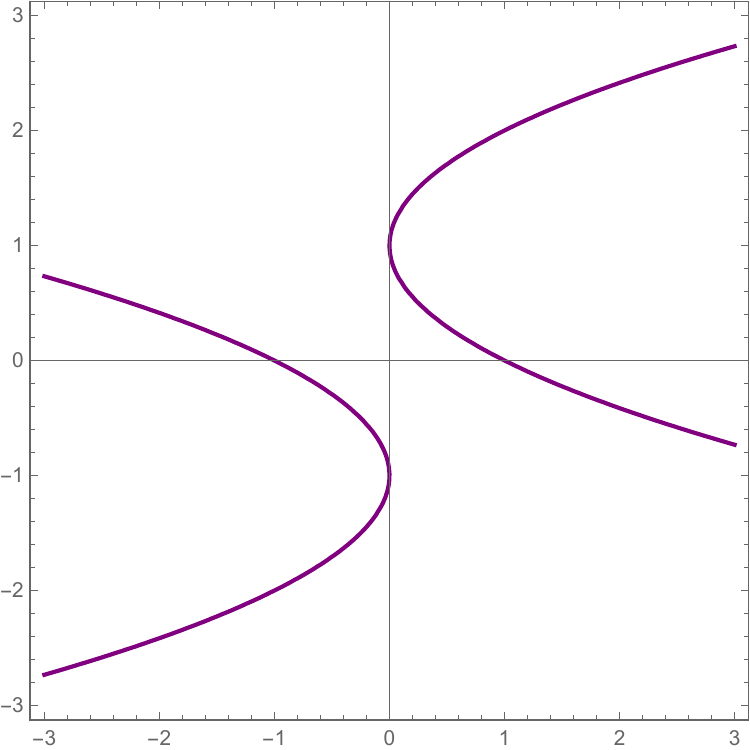}
			\caption{Without constant geometric fiber cardinality} 
			\label{fig:real-fiber-insufficient}
		\end{subfigure}
		~
		\caption{Non-examples}
	\end{figure}	
	\begin{enumerate}[label=(\alph*)]
		\item The flatness assumption on the morphism is essential and natural. Consider the variety defined by the union of the hyperbola $V(xy-1)$ and the origin $V(x,y)$ on the plane and its projection to the affine line. Obviously there is no section to the projection near the origin. Notice that the geometric fiber size is always 1. In this example, the projection is not flat, since it is not open. See Figure \ref{fig:non-flatness}.
		
		\item Also, the constancy assumption on the geometric fiber cardinality is crucial. Let $X=V((y-1)^2-x)\cup V((y+1)^2+x)$ be the union of two parabolas and let $Y=\spec R[x]$ be the affine line, then the canonical projection from $X$ to $Y$ is obviously flat. The set $X(R)$ is visualized in Figure \ref{fig:real-fiber-insufficient}. Clearly there are no two distinct semi-algebraic continuous sections near the origin.	Notice that the number of real fibers over $Y(R)$ is always $2$ in this example. So even a constancy condition on the real fiber is not enough.
		
	\end{enumerate}		
\end{example}

The complex analogue of Theorem \ref{thm-covering-map} is also true. 

\begin{theorem}
	Let $X$, $Y$ be two $C$-varieties. Suppose $\pi: X\to Y$ is q-\'etale $C$-morphism. Then $\pi_C$ is a covering map in the Euclidean topology. %
\end{theorem}
\begin{proof}
	Replacing $\pi:X\to Y$ with the reduced map, we may assume that $\pi$ is finite and \'etale. Because the question is local, we may also assume that $X$, $Y$ are affine varieties. Since $C/R$ is a finite field extension, the Weil restriction $\mathop{\mathrm{Res}_{C/R}}$ exists for quasi-projective varieties. Let $Z=\mathop{\mathrm{Res}_{C/R}}{X}$ and $W=\mathop{\mathrm{Res}_{C/R}}{Y}$ be the Weil restrictions of $X$ and $Y$ respectively, then $Z$ and $W$ are $R$-varieties whose $R$-rational points can be naturally identified with closed points of $X$ and $Y$. Also, the induced map $\mathrm{Res}_{C/R}(\pi):Z\to W$ is finite and \'etale \cite[Proposition 4.9, Proposition 4.10.1]{scheiderer2006real}. The proof is completed by applying Theorem \ref{thm-covering-map} to $\mathrm{Res}_{C/R}(\pi):Z\to W$.
\end{proof}

\begin{remark}
	One possible way to understand Theorem \ref{thm-covering-map} is to compare it with the Ehresmann's fibration theorem \cite{ehresmann1950connexions}. Actually,
	Theorem \ref{thm-covering-map} implies a weak semi-algebraic variant of Ehresmann's fibration theorem. Suppose that $\pi:X\to Y$ is a proper morphism of non-singular $R$-varieties with surjective tangent maps $T_\pi: T_x X\to T_{\pi(x)} Y \otimes_{\kappa_y} \kappa_x$, and suppose further that all fibers of $\pi$ are finite, then $\pi_R$ is a covering map in the Euclidean topology. Indeed, such a morphism is  quasi-finite, proper and \'etale \cite[Corollary VII.5.5]{altman2006introduction}. Therefore locally there are semi-algebraic sections to $\pi_R$ by Theorem \ref{thm-covering-map}.
\end{remark}

\section{Applications}\label{sect-app}

We conclude the paper by showing some consequences of Theorem \ref{thm-covering-map}.

In \cite{hardt1980semi}, Hardt proved that a semi-algebraic continuous map $f$ between two semi-algebraic sets $X,Y$ can be stratified so that the map is locally trivial on each stratum, which is called Hardt's Triviality. His strategy is to construct a $(\dim Y-1)$-dimensional semi-algebraic subset $Z$ of $Y$ so that $f$ is locally trivial outside $Z$. The proof factorizes $f$ as consecutive corank $1$ projections of Euclidean space to deform the fibers, so $Z$ might be necessarily enlarged to control some ramification that occurs in the middle only. He then asked if it is possible to construct a smallest $Z$. Here we present a different triviality result for a quasi-finite morphism without using successive corank 1 projections, so it is more likely to generate a smaller $Z$.

\begin{proposition}\label{cor-triviality}
	Suppose $f:X\to Y$ is a quasi-finite $R$-morphism of $R$-varieties, then there are finitely many locally closed subschemes $\{Y_i\}$ of $Y$, whose union covers $Y$, such that the restriction map $f_i:X_i=X\times_Y Y_i\to Y_i$ is a covering map on the $R$-rational points.
\end{proposition}
\begin{proof}
	By taking reduced subscheme structures and applying generic flatness repeatedly, there exists $t\geq 0$ and closed subschemes $$Y\supseteq Y_{\mathrm{red}}=V_0\supsetneq V_1\supsetneq\cdots\supsetneq V_t=\varnothing$$ such that $X\times_Y (V_i\backslash V_{i+1})$ is flat over $V_i\backslash V_{i+1}$ (the generic flatness stratification, see \cite[\href{https://stacks.math.columbia.edu/tag/0H3Z}{Lemma 0H3Z}]{stacks-project}). Set $Y_i=V_i\backslash V_{i+1}$, then the restriction of $f$ on $Y_i$ is quasi-finite and flat. By the lower semi-continuity of geometric fiber cardinality (Proposition \ref{prop-grothendieck}), we can further stratify $Y_i$ into locally closed subschemes $Y_{ij}$ where the cardinality of geometric fibers is a constant. Hence the restriction of $f$ is q-\'etale on $Y_{ij}$. Now apply Theorem \ref{thm-covering-map} and the proof is completed.
\end{proof}

\begin{example}
	Set $X=\spec \mathbb{R}[x,y,z,w]/\ideal{x^2+y^2-1,z^2+w^2-1,x-2z^2+1,x+2 w^2-1,y-2 z w}$ and $Y=\spec \mathbb{R}[x,y]/\ideal{x^2+y^2-1}$. In Example \ref{examples-of-covering-maps}, we saw that $\pi:X\to Y$ induces a covering map on the $\mathbb{R}$-rational points, because $\pi$ itself is q-\'etale. 
	
	If $\pi$ is factored as successive corank $1$ projections $X\to\mathbb{A}_\mathbb{R}^3\to\mathbb{A}_\mathbb{R}^2$, then some ramification occurs in the middle. Denote the middle image by $Z=\spec \mathbb{R}[x,y,z]/\ideal{x^2+y^2-1,x-2z^2+1}$, then $\pi$ is the composition of $\pi_1:X\to Z$ and $\pi_2:Z\to Y$. Generically $\pi_2$ is a 2-to-1 map, but there is a ramification over $(-1,0)$. As a result, $(-1,0,0,\pm 1)\mapsto (-1,0,0)\mapsto (-1,0)$ has to be treated separately in the composition $\pi=\pi_2\circ \pi_1$. One way to see this is to think $\pi_\mathbb{R}$ as the double cover of $S^1$, and think $Z$ as an $\infty$-shaped curve.
	
\end{example}

Now we provide an explanation of the core mechanism of cylindrical algebraic decomposition. In cylindrical algebraic decomposition, the delineability of roots of a polynomial $f\in R[x_1,\ldots,x_{r-1}][x_r]$ can be concluded from the sign-invariance of the coefficients and sub-discriminants. To be more precisely, we restate \cite[Theorem 4]{collins1975quantifier} here.

\begin{theorem}[Collins]
	Let $A(x_1,\ldots,x_r)=\sum_{k=0}^{n} c_k(x_1,\ldots,x_{r-1})x_r^k$ be a real polynomial, $r\geq 2$, $S$ a (semi-algebraically) connected subset of $R^{r-1}$. Let $$\mathcal{B}=\left\{\sum_{k=0}^{j} c_k(x_1,\ldots,x_{r-1})x_r^k\middle|0\leq j\leq n\right\}.$$ and $$\mathcal{P}=\left\{c_k(x_1,\ldots,x_{r-1})\middle|0\leq k\leq n\right\}\cup \left\{\mathtt{sRes}_{x_r,j}(B,\frac{\partial B}{\partial x_r})\middle|B\in \mathcal{B}, 0\leq j <\deg_{x_r} B\right\}.$$ If every element of $\mathcal{P}$ is sign-invariant on $S$, then the roots of $A$ are delineable on $S$. That is, the real roots of $A(x_1,\ldots,x_{r-1};x_r)=0$ are continuous semi-algebraic functions in $(x_1,\ldots,x_{r-1})$.
\end{theorem}

Let $V=V(A)$ be the hypersurface defined by $A$ in $\mathbb{A}_R^r$. Consider the projection $\pi:V(A)\to \mathbb{A}_R^{r-1}$, then a generic flatness stratification of $\pi$ is given by

$$\mathbb{A}_R^{r-1}\supseteq V(c_n) \supseteq V(c_n,c_{n-1})\supseteq V(c_n, c_{n-1},c_{n-2})\supseteq \cdots \supseteq V(c_n,\ldots,c_0)\supseteq \varnothing.$$

Let $Y_0= \mathbb{A}_R^{r-1}\backslash V(c_n)$, $Y_1=V(c_n) \backslash V(c_n,c_{n-1})$, $\ldots$, $Y_n= V(c_n,\ldots,c_1) \backslash V(c_n,\ldots,c_0)$, $Y_{n+1}=V(c_n,\ldots,c_0)$.
Over each stratum $Y_i$, the projection is flat. Moreover, the projection is finite over every strata except the last one $Y_{n+1}$. This explains the coefficients $c_i$ in Collins's projection $\mathcal{P}$. Then, in the proof of Proposition \ref{cor-triviality}, $Y_0,\ldots,Y_n$ are further partitioned by the number of geometric points in the fiber. This number is detected by the sub-discriminants $\mathtt{sRes}_{x_r,j}(B,\frac{\partial B}{\partial x_r})$ in $\mathcal{P}$. Though the projection is not finite over $Y_{n+1}$, no extra care is needed because over $Y_{n+1}$, $V(A)$ is simply a cylinder. 

Therefore there are sufficient polynomials in $\mathcal{P}$ to encode a q-\'etale stratification of $\pi$. By Proposition \ref{cor-triviality}, $\pi$ induces a covering on $S$. Because the sections can be ordered in $x_r$, local semi-algebraic sections of $\pi$ on $S$ can be extended to whole $S$. Now we arrive at the conclusion that the real roots of $A$ are delineable over $S$.

Let us turn to a related problem. Finding at least one point on each semi-algebraically connected component of a given real algebraic set is a fundamental problem in computational real algebraic geometry. Usually, the most practical method to solve this is computing a cylindrical algebraic decomposition, which involves successive corank 1 projections. We point out that one can reduce the problem to several smaller similar problems, using the triviality result we just established. %

\begin{proposition}
	Suppose $f:X\to Y$ is a quasi-finite $R$-morphism of $R$-varieties, then finding at least one sample point on each semi-algebraically connected component of $X(R)$ can be reduced to finding at least one sample point on each semi-algebraically connected component of $Y_i(R)$ for some locally closed subschemes $Y_i\subseteq Y$.
\end{proposition}
\begin{proof}
	By Corollary \ref{cor-triviality}, there are finitely many locally closed subschemes $Y_i\subseteq Y$ such that $X_i=X \times_Y Y_i \to Y_i$ induces a covering map on the rational points. Clearly, if we have found at least one sample point on each semi-algebraically connected component of each $X_i(R)$, then the union of sample points meets every semi-algebraically connected component of $X(R)$. So it suffices to show that sample points meeting each semi-algebraically connected components of $Y_i(R)$ can be lifted to sample points meeting each semi-algebraically connected components of $X_i(R)$. Let $S=\{y_i\}\subseteq Y_i(R)$ be the sample points of $Y_i(R)$, we claim that $f^{-1}_R(S)=\bigcup_i f_R^{-1}(y_i)$ is a set of sample points meeting each semi-algebraically connected components of $X_i(R)$.
	
	Let $C\subseteq X_i(R)$ be a semi-algebraically connected component of $X_i(R)$. Pick arbitrary $x_0\in C$ and let $y_0=f(x_0)$. Then there is some $y_i\in S$ such that $y_0$ and $y_i$ lie in the same semi-algebraically connected component. Notice that for a semi-algebraic set, being semi-algebraically connected is equivalent to being semi-algebraically path-connected \cite[Proposition 2.5.13]{bochnak2013real}. So there is a semi-algebraic path $\gamma:[0,1]\to Y_i(R)$ such that $\gamma(0)=y_0$ and $\gamma(1)=y_i$. By the unique path lifting property \cite[Proposition 5.7]{delfs1984introduction}, there is a semi-algebraic path $\tilde{\gamma}:[0,1] \to X_i(R)$ such that $\tilde{\gamma}(0)=x_0$ and $f_R\circ \tilde{\gamma}=\gamma$. Then $\tilde{\gamma}(1)\in f^{-1}_R(y_i)$ is the required sample point meeting $C$.
\end{proof}

\printbibliography
\end{document}